\documentclass[11pt]{article}

\usepackage{amsbsy}
\usepackage{amssymb}
\usepackage{amsthm}
\usepackage{amsmath}
\usepackage{fullpage,hyperref,mathrsfs,txfonts,epsfig,graphicx,graphics,verbatim}

\begin{document}

\newcommand{\bwr}{\boldsymbol{\wr}}

% ENVIRONMENTS & THEOREMS

\newenvironment{nmath}{\begin{center}\begin{math}}{\end{math}\end{center}}

\newtheorem{thm}{Theorem}[section]
\newtheorem{lem}[thm]{Lemma}
\newtheorem{remark}[thm]{Remark}
\newtheorem{prop}[thm]{Proposition}
\newtheorem{cor}[thm]{Corollary}
\newtheorem{conj}[thm]{Conjecture}
\newtheorem{dfn}[thm]{Definition}
\newtheorem{prob}[thm]{Problem}
\newtheorem{ques}[thm]{Question}

% COMMANDS

\newcommand{\A}{\mathcal{A}}
\newcommand{\B}{\mathcal{B}}
\newcommand{\K}{\mathcal{K}}
\newcommand{\F}{\mathbb{F}}
\newcommand{\p}{\Psi}
\newcommand{\f}{\Phi}
\newcommand{\Lip}{\mathrm{Lip}}
\newcommand{\1}{\mathbf{1}}
\newcommand{\s}{\sigma}
\renewcommand{\P}{\mathcal{P}}
\renewcommand{\O}{\Omega}
\renewcommand{\S}{\Sigma}
\renewcommand{\Pr}{\mathbb{P}}
\renewcommand{\approx}{\asymp}
\newcommand{\T}{\mathrm{T}}
\newcommand{\co}{\mathrm{co}}
\newcommand{\Isom}{\mathrm{Isom}}
\newcommand{\edge}{\mathrm{edge}}
\newcommand{\bounded}{\mathrm{bounded}}
\newcommand{\e}{\varepsilon}
\newcommand{\im}{\mathrm{i}}
\newcommand{\restrict}{\upharpoonright}
\newcommand{\supp}{{\mathrm{\bf supp}}}
\renewcommand{\l}{\lambda}
\newcommand{\U}{\mathcal{U}}
\newcommand{\calH}{\mathcal{H}}
\newcommand{\G}{\Gamma}
\newcommand{\g}{\gamma}
\renewcommand{\L}{\mathscr{L}}
\newcommand{\hcf}{\mathrm{hcf}}
\renewcommand{\a}{\alpha}
\newcommand{\N}{\mathbb{N}}
\newcommand{\R}{\mathbb{R}}
\newcommand{\Z}{\mathbb{Z}}
\newcommand{\C}{\mathbb{C}}

\newcommand{\E}{\mathbb{E}}
\newcommand{\alp}{\alpha^*}

\newcommand{\bb}[1]{\mathbb{#1}}
\renewcommand{\rm}[1]{\mathrm{#1}}
\renewcommand{\cal}[1]{\mathcal{#1}}

\newcommand{\fin}{\nolinebreak\hspace{\stretch{1}}$\lhd$}

\title{The Johnson-Lindenstrauss lemma almost characterizes\\ Hilbert space, but not quite}
\author{William B. Johnson \footnote{Research
supported in part by NSF grants DMS-0503688 and DMS-0528358.}\\Texas A\& M
University\\{\tt johnson@math.tamu.edu}\and Assaf
Naor\footnote{Research supported in part by NSF grants CCF-0635078
and DMS-0528387.}\\Courant Institute\\{\tt naor@cims.nyu.edu} }
\date{}

\maketitle

\begin{abstract}
Let $X$ be a normed
%%%%CHANGE
space that satisfies the Johnson-Lindenstrauss
lemma
%%%%CHANGE
(J-L lemma, in short)
in the sense that for any integer $n$ and any $x_1,\ldots,x_n\in X$ there
exists a linear mapping $L:X\to F$, where $F\subseteq X$ is a linear
subspace of dimension $O(\log n)$, such that $\|x_i-x_j\|\le
\|L(x_i)-L(x_j)\|\le O(1)\cdot\|x_i-x_j\|$ for all $i,j\in
\{1,\ldots, n\}$. We show that this implies that $X$ is almost
Euclidean in the following sense:
%%%%CHANGE
Every $n$-dimensional subspace of
$X$ embeds into Hilbert space with distortion $2^{2^{O(\log^*n)}}$.
On the other hand, we show that there exists a normed space $Y$
which satisfies the
%%%%CHANGE
J-L lemma, but for every $n$
there exists an $n$-dimensional subspace $E_n\subseteq Y$ whose
Euclidean distortion is at least $2^{\Omega(\alpha(n))}$, where
$\alpha$ is the inverse Ackermann function.
\end{abstract}

\section{Introduction}

The
%%%%CHANGE
J-L lemma~\cite{JL84} asserts that if $H$ is a
Hilbert space, $\e>0$, $n\in \N$, and $x_1,\ldots,x_n\in H$ then
there exists a linear mapping (even a multiple of an orthogonal
projection) $L:H\to F$, where $F\subseteq H$ is a linear subspace of
dimension $O(c(\e)\log n)$, such that for all $i,j\in
\{1,\ldots,n\}$ we have \begin{eqnarray}\label{eq:JL}\|x_i-x_j\|\le
\|L(x_i)-L(x_j)\|\le (1+\e)\|x_i-x_j\|.
\end{eqnarray}
This fact has found many applications in mathematics and computer
science, in addition to the original application in~\cite{JL84} to
%%%%CHANGE
a Lipschitz extension problem. The widespread applicability of the
%%%%CHANGE
J-L lemma in computer science can be (somewhat
simplistically) attributed to the fact that it can be viewed as a
compression scheme which helps
%%%%CHANGE
to reduce significantly  the space
required for storing multidimensional data. We shall not attempt to
list here all the applications of the
%%%%CHANGE
J-L lemma to
areas ranging from nearest neighbor search to machine learning---we
refer the interested reader to~\cite{Kleinberg97,IM99,
KOR00,Indyk01,Vem04,Indyk06,AC06} and the references therein for a
partial list of such applications.

The applications of~\eqref{eq:JL} involve various requirements from
the mapping $L$. While some applications just need the distance
preservation condition~\eqref{eq:JL} and not the linearity of $L$,
most
%%%%CHANGE
applications require $L$ to be linear. Also, many applications
are based on additional information that comes from the proof of the
%%%%CHANGE
J-L  lemma, such as the fact that $L$ arises with
high probability from certain distributions over linear mappings.
The linearity of $L$ is useful, for example, for fast evaluation of
the images $L(x_i)$, and also because these images behave well when
additive noise is applied to the initial vectors $x_1,\ldots,x_n$.

Due to the usefulness of the
%%%%CHANGE
J-L lemma there has
been considerable effort by researchers to prove such a
dimensionality reduction theorem in other normed spaces. All of
these efforts have thus far resulted in negative results which show
that the
%%%%CHANGE
J-L lemma fails to hold true in certain non-Hilbertian settings.
In~\cite{CS02} Charikar and Sahai proved that there is no dimension
reduction via linear mappings in $L_1$. This negative result was
extended to any $L_p$, $p\in [1,\infty]\setminus\{2\}$, by Lee,
Mendel and Naor in~\cite{LMN05}. Negative results for dimension
reduction without the requirement that the embedding $L$ is linear
are known only for the spaces $L_1$~\cite{BC05,LN04,LMN05}   and
$L_\infty$~\cite{Bou85, JLS87, AR92, Mat96, LMN05}.
%%%%CHANGE see~\cite{AR92,Mat96,BC05,LN04,LMN05}.
Here we show
that the negative results for linear dimension reduction in $L_p$
spaces are a particular case of a much more general phenomenon: A
normed space that satisfies the
%%%%CHANGE
J-L lemma is very
close to being Euclidean in the sense that all of its
$n$-dimensional subspaces are isomorphic to Hilbert space with
distortion $2^{2^{O(\log^*(n))}}$. Here, and in what follows, if
$x\ge 1$ then $\log^*(x)$ is the unique integer $k$ such that if we
define $a_{1}=1$ and $a_{i+1}=e^{a_i}$ (i.e. $a_i$ is an exponential
tower of height $i$), then $a_{k}< x\le a_{k+1}$.

In order to state our results we recall the following notation: The
Euclidean distortion of a finite dimensional normed space $X$,
denoted $c_2(X)$, is the infimum over all $D>0$ such that there
exists a linear mapping $S:X\to \ell_2$ which satisfies $\|x\|\le
\|S(x)\|\le D\|x\|$ for all $x\in X$. Note that in the computer
science literature the notation $c_2(X)$ deals with bi-Lipschitz
embeddings, but in the context of normed spaces it can be shown that
the optimal bi-Lipschitz embedding may be chosen to be linear (this
is explained for example
%%%%CHANGE
 in~\cite[Chapter 7]{BL00}). The parameter $c_2(X)$ is
also known as the Banach-Mazur distance between $X$ and Hilbert
space.

\begin{thm}\label{thm:upper} For every $D,K>0$ there exists a
constant $c=c(K,D)>0$ with the following property. Let $X$ be a
Banach space such that for every $n\in \N$ and every
$x_1,\ldots,x_n\in X$ there exists a linear subspace $F\subseteq X$,
of dimension at most $K\log n$, and a linear mapping $S:X\to F$ such
that $\|x_i-x_j\|\le \|S(x_i)-S(x_j)\|\le D\|x_i-x_j\|$ for all
$i,j\in \{1,\ldots,n\}$. Then for every $k\in \N$ and every
$k$-dimensional subspace $E\subseteq X$, we have
\begin{eqnarray}\label{eq:log*}
c_2(E)\le 2^{2^{c\log^*(k)}}.
\end{eqnarray}
\end{thm}

The proof of Theorem~\ref{thm:upper} builds on ideas
from~\cite{CS02,LMN05}, while using several fundamental results from
the local theory of Banach spaces. Namely, in~\cite{LMN05} the $L_1$
point-set from~\cite{CS02} was analyzed via an analytic argument
which extends to any $L_p$ space, $p\neq 2$, rather than the linear
programming argument in~\cite{CS02}. In Section~\ref{sec:upper} we
construct a variant of this point-set in any Banach space, and use
it in conjunction with some classical results in Banach space theory
to prove Theorem~\ref{thm:upper}.

\medskip

%%%%CHANGE
The fact that the bound on
 $c_2(E)$  in~\eqref{eq:log*} is not
$O(1)$ is not just an artifact of our iterative proof technique:
There do exist non-Hilbertian Banach spaces which satisfy the
%%%%CHANGE
J-L lemma!

\begin{thm}\label{thm:lower}
There exist two universal constants $D,K>0$ and a Banach space $X$
such that for every $n\in \N$ and every $x_1,\ldots,x_n\in X$ there
exists a linear subspace $F\subseteq X$, of dimension at most $K\log
n$, and a linear mapping $S:X\to F$ such that $\|x_i-x_j\|\le
\|S(x_i)-S(x_j)\|\le D\|x_i-x_j\|$ for all $i,j\in \{1,\ldots,n\}$.
Moreover, for every integer $n$ the space $X$ has an $n$-dimensional
subspace $F_n\subseteq X$ with
\begin{eqnarray}\label{eq:bellenot}
c_2(F_n)\ge 2^{c\alpha(n)},
\end{eqnarray}
where $c>0$ is a universal constant and $\alpha(n)\to\infty$ is the
inverse Ackermann function.
\end{thm}
We refer the readers to Section~\ref{sec:lower} for the definition
of the inverse Ackermann function. The Banach space $X$ in
Theorem~\ref{thm:lower} is the $2$-convexification of the Tsirelson
space~\cite{Ts74}, denoted $T^{(2)}$, which we shall now define.
%%%%CHANGE
The definition below, due to Figiel and Johnson~\cite{FJ74}, actually gives the dual to the
space constructed by Tsirelson (see the
book~\cite{CS89} for a comprehensive discussion). Let $c_{00}$ denote
the space of all finitely supported sequences of real numbers. The
standard unit basis of $c_{00}$ is denoted by
$\{e_i\}_{i=1}^\infty$. Given $A\subseteq \N$ we denote by $P_A$ the
restriction operator to $A$, i.e. $P_A\left(\sum_{i=1}^\infty
x_ie_i\right)=\sum_{i\in A} x_i e_i$. Given two finite subsets
$A,B\subseteq \N$ we write $A<B$
%(resp. $A\le B$)
if $\max A<\min B$
%(resp. $\max A\le \min B$)
.
Define inductively a sequence of norms
$\{\|\cdot\|_m\}_{m=0}^\infty$ by $\|x\|_0=\|x\|_{c_0}=\max_{j\ge 1}
|x_j|$, and
\begin{eqnarray}\label{eq:tsirl recurse}
\|x\|_{m+1}=\max\left\{
\|x\|_{m},\frac12\sup\left\{\sum_{j=1}^n\left\|P_{A_j}(x)\right\|_m:\
n\in \N,\ A_1,\ldots,A_n\subseteq \N\ \mathrm{finite,}\  \{n\}<
A_1<A_2<\cdots<A_n\right\}\right\}
\end{eqnarray}
Then for each $x\in c_{00}$ the sequence $\{\|x\|_m\}_{m=0}^\infty$
is nondecreasing and bounded from above by
$\|x\|_{\ell_1}=\sum_{j=1}^\infty |x_j|$. It follows that the limit
$\|x\|_{T}\coloneqq \lim_{m\to\infty} \|x\|_m$ exists. The space
$X=T^{(2)}$ from Theorem~\ref{thm:lower} is the completion of
$c_{00}$ under the norm:
\begin{eqnarray}\label{eq:T2}
\left\|\sum_{j=1}^\infty x_j e_j\right\|_{T^{(2)}}\coloneqq
\left\|\sum_{j=1}^\infty |x_j|^2 e_j\right\|_{T}^{1/2}.
\end{eqnarray}
The proof of the fact that $T^{(2)}$ satisfies the
%%%%CHANGE
J-L lemma consists of a concatenation of several
classical results, some of which are quite deep. The lower
bound~\eqref{eq:bellenot} follows from the work of
Bellenot~\cite{Bell84}. The details are presented in
Section~\ref{sec:lower}.

\section{Proof of Theorem~\ref{thm:upper}}\label{sec:upper}

Let $(X,\|\cdot\|)$ be a normed space. The Gaussian type $2$
constant of $X$, denoted $T_2(X)$, is the infimum over all $T>0$
such that for every $n\in \N$ and every $x_1,\ldots,x_n\in X$ we
have,
\begin{eqnarray}\label{eq:def type}
\E\left\|\sum_{i=1}^n g_ix_i\right\|^2\le T^2 \sum_{i=1}^n
\|x_i\|^2.
\end{eqnarray}
Here, and in what follows, $g_1,\ldots,g_n$ denote i.i.d. standard
Gaussian random variables. The cotype $2$ constant of $X$, denoted
$C_2(X)$, is the infimum over all $C>0$ such that for every $n\in
\N$ and every $x_1,\ldots,x_n\in X$ we have,
\begin{eqnarray}\label{eq:def cotype}
\sum_{i=1}^n \|x_i\|^2\le C^2 \E\left\|\sum_{i=1}^n
g_ix_i\right\|^2.
\end{eqnarray}

A famous theorem of Kwapien~\cite{Kw72} (see also the exposition
%%%%CHANGE
in~\cite[Theorem 3.3]{Pisier86}) states that
\begin{eqnarray}\label{eq:kwapien}
c_2(X)\le T_2(X)\cdot C_2(X).
\end{eqnarray}
An important theorem of Tomczak-Jaegermann~\cite{T-J79} states that
if the Banach space $X$ is $d$-dimensional then there exist
$x_1,\ldots,x_d,y_1\ldots,y_d\in X\setminus\{0\}$ for which
\begin{eqnarray}\label{eq:T-J}
\E\left\|\sum_{i=1}^d g_ix_i\right\|^2\ge \frac{T_2(X)^2}{2\pi}
\sum_{i=1}^d \|x_i\|^2\quad \mathrm{and}\quad \sum_{i=1}^d
\|y_i\|^2\ge \frac{C_2(X)^2}{2\pi} \E\left\|\sum_{i=1}^d
g_iy_i\right\|^2.
\end{eqnarray}
In other words, for $d$-dimensional spaces it suffices to consider
$n=d$ in~\eqref{eq:def type} and~\eqref{eq:def cotype} in order to
compute $T_2(X)$ and $C_2(X)$ up to a universal factor. For our
purposes it suffices to use the following simpler fact due to
Figiel, Lindenstrauss and
%%%%CHANGE
Milman~\cite[Lemma 6.1]{FLM77}:
If $\dim(X)=d$ then there exist
$x_1,\ldots,x_{d(d+1)/2},y_1,\ldots,y_{d(d+1)/2}\in X\setminus\{0\}$
for which
\begin{eqnarray}\label{eq:FLM}
\E\left\|\sum_{i=1}^{d(d+1)/2} g_ix_i\right\|^2= T_2(X)^2
\sum_{i=1}^{d(d+1)/2} \|x_i\|^2\quad \mathrm{and}\quad
\sum_{i=1}^{d(d+1)/2} \|y_i\|^2= C_2(X)^2
\E\left\|\sum_{i=1}^{d(d+1)/2} g_iy_i\right\|^2.
\end{eqnarray}
We note, however, that it is possible to improve the constant terms
in Theorem~\ref{thm:upper} if we use~\eqref{eq:T-J} instead
of~\eqref{eq:FLM} in the proof below. We shall now sketch the proof
of~\eqref{eq:FLM}, taken
%%%%CHANGE
from ~\cite[Lemma 6.1]{FLM77}, since this
type of finiteness result is used crucially in our proof of
Theorem~\ref{thm:upper}.

We claim that if $m>d(d+1)/2$ and $u_1,\ldots,u_m\in X$ then there
are $v_1,\ldots,v_{m-1}, w_1,\ldots,w_{m-1}\in X$ such that
\begin{eqnarray}\label{eq:decomposition}
\E\left\|\sum_{i=1}^{m-1} g_iv_i\right\|^2+\E\left\|\sum_{i=1}^{m-1}
g_iw_i\right\|^2=\E\left\|\sum_{i=1}^{m} g_iu_i\right\|^2 \quad
\mathrm{and}\quad \sum_{i=1}^{m-1} \|v_i\|^2+\sum_{i=1}^{m-1}
\|w_i\|^2=\sum_{i=1}^m \|u_i\|^2.
\end{eqnarray}
Note that~\eqref{eq:decomposition} clearly implies~\eqref{eq:FLM}
since it shows that in the definitions~\eqref{eq:def type}
and~\eqref{eq:def cotype} we can take $n=d(d+1)/2$ (in which case
the infima in these definitions are attained by a simple compactness
argument).

To prove~\eqref{eq:decomposition} we can think of $X$ as $\R^d$,
equipped with a norm $\|\cdot\|$. The random vector $\sum_{i=1}^m
g_i u_i=\left(\sum_{i=1}^d g_iu_{ij}\right)_{j=1}^d$ has a Gaussian
distribution with covariance matrix $A=(\sum_{i=1}^m
u_{ij}u_{ik})_{j,k=1}^d=\sum_{i=1}^m u_i\otimes u_i$. Thus the
symmetric matrix $A$ is in the cone generated by the symmetric
matrices $\{u_i\otimes u_i\}_{i=1}^m$. By Caratheodory's theorem for
cones (see e.g. \cite{DGK63}) we may reorder the vectors $u_i$ so as
to find scalars $c_1\ge c_2\ge\cdots\ge c_m\ge 0$ with $c_i=0$ for
$i>d(d+1)/2$, such that $A=\sum_{i=1}^m c_i u_i\otimes u_i$. This
sum contains at most $d(d+1)/2$ nonzero summands. Define
$v_i\coloneqq \sqrt{c_i/c_1}\cdot u_i$ (so that there are at most
$d(d+1)/2\le m-1$ nonzero $v_i$) and $w_i=\sqrt{1-c_i/c_1}\cdot u_i$
(so that $w_1=0$). The second identity in~\eqref{eq:decomposition}
is trivial with these definitions. Now, the random vector
$\sum_{i=1}^m g_iv_i$ has covariance matrix
$\frac{1}{c_1}\sum_{i=1}^m c_iu_i\otimes u_i=\frac{1}{c_1} A$ and
the random vector $\sum_{i=1}^m g_iw_i$ has covariance matrix
$\sum_{i=1}^m (1-c_i/c_1)u_i\otimes u_i=(1-1/c_1)A$. Thus
$\E\left\|\sum_{i=1}^{m-1}
g_iv_i\right\|^2=\frac{1}{c_1}\E\left\|\sum_{i=1}^{m}
g_iu_i\right\|^2$ and $\E\left\|\sum_{i=1}^{m-1}
g_iw_i\right\|^2=(1-1/c_1)\E\left\|\sum_{i=1}^{m} g_iu_i\right\|^2$.
This completes the proof of~\eqref{eq:decomposition}.

\medskip

We are now in position to prove Theorem~\ref{thm:upper}. Define
\begin{eqnarray}\label{eq:def Delta}
\Delta(n)\coloneqq \Delta_X(n)\coloneqq \sup\left\{c_2(F):\ F\subseteq X\ \mathrm{linear\
subspace,}\ \dim(F)\le n\right\}.
\end{eqnarray}
Note that by John's theorem~\cite{John48} (see also the beautiful
exposition in~\cite{Ball97}) $\Delta(n)\le \sqrt{n}$. Our goal is to
obtain a much better bound on $\Delta(n)$. To this end let
$F\subseteq X$ be a linear subspace of dimension $k\le n$. Let $m$
be the integer satisfying $2^{m-1}<k(k+1)/2\le 2^m$. We shall use
the vectors from~\eqref{eq:FLM}. By adding some zero vectors so as
to have $2^m$ vectors, and labeling them (for convenience) by the
subsets of $\{1,\ldots,m\}$, we obtain $\{x_A\}_{A\subseteq
\{1,\ldots,m\}},\{y_A\}_{A\subseteq \{1,\ldots,m\}}\subseteq X$ such
that
%%%%CHANGE X to F
\begin{eqnarray}\label{eq:relabeledx}
\E\left\|\sum_{A\subseteq \{1,\ldots,m\}} g_Ax_A\right\|^2= T_2(F)^2
\sum_{A\subseteq \{1,\ldots,m\}} \|x_A\|^2>0 \end{eqnarray}
\begin{eqnarray}\label{eq:relabeledy}
\sum_{A\subseteq \{1,\ldots,m\}} \|y_A\|^2= C_2(F)^2
\E\left\|\sum_{A\subseteq \{1,\ldots,m\}} g_Ay_A\right\|^2>0.
\end{eqnarray}
For every $\e=(\e_1,\ldots,\e_m)\in \{-1,1\}$ and $A\subseteq
\{1,\ldots,m\}$ consider the Walsh function $W_A(\e)=\prod_{i\in A}
\e_i$. For every $g=\{g_A\}_{A\subseteq \{1,\ldots,m\}}$ define
$\f_g,\p_g:\{-1,1\}^m\to F$ by
\begin{eqnarray}\label{eq:def functions}
\f_g(\e)\coloneqq \sum_{A\subseteq \{1,\ldots,m\}} g_A W_A(\e)
x_A\quad \mathrm{and}\quad \p_g(\e)\coloneqq \sum_{A\subseteq
\{1,\ldots,m\}} g_A W_A(\e) y_A.
\end{eqnarray}
Thus $\f_g,\p_g$ are random $F$-valued functions given by the random
Fourier expansions in~\eqref{eq:def functions}---the randomness is
with respect to the i.i.d. Gaussians $g=\{g_A\}_{A\subseteq
\{1,\ldots,m\}}$. These random functions induce the following two
random subsets of $F$:
$$
U_g\coloneqq \{\f_g(\e)\}_{\e\in \{-1,1\}^m}\cup \{x_A\}_{A\subseteq
\{1,\ldots,m\}}\cup \{0\}\quad\mathrm{and}\quad V_g\coloneqq
\{\p_g(\e)\}_{\e\in \{-1,1\}^m}\cup \{y_A\}_{A\subseteq
\{1,\ldots,m\}}\cup \{0\}.
$$
Then $|U_g|,|V_g|\le 2^{m+1}+1\le 2k(k+1)+1\le 2(n+1)^2$. By the
assumptions of Theorem~\ref{thm:upper} it follows that there exist
two subspaces $E_g,E_g'\subseteq X$ with $\dim(E_g),\dim(E_g')\le
K\log\left(2(n+1)^2\right)\le 4K\log(n+1)$ and two linear mappings
$L_g:X\to E_g$, $L_g':X\to E_g'$, which satisfy
\begin{eqnarray}\label{eq:bilip}
x,y\in U_g\implies \|x-y\|\le \|L_g(x)-L_g(y)\|\le D\|x-y\|.
\end{eqnarray}
and
\begin{eqnarray}\label{eq:bilip'}
x,y\in V_g\implies \|x-y\|\le \|L_g'(x)-L_g'(y)\|\le D\|x-y\|.
\end{eqnarray}
Moreover, by the definition of $\Delta(\cdot)$ there are two linear
mappings $S_g:E_g\to \ell_2$ and $S_g':E_g'\to \ell_2$ which satisfy
\begin{eqnarray}\label{eq:euclidean}
x\in E_g\implies \|x\|\le \|S_g(x)\|_2\le
2\Delta\left(4K\log(n+1)\right)\|x\|,
\end{eqnarray}
and
\begin{eqnarray}\label{eq:euclidean'}
x\in E_g'\implies \|x\|\le \|S_g'(x)\|_2\le
2\Delta\left(4K\log(n+1)\right)\|x\|.
\end{eqnarray}
By the orthogonality of the Walsh functions we see that
\begin{eqnarray}\label{eq:orthogonal}
\E_\e
\left\|S_g\left(L_g\left(\f_g(\e)\right)\right)\right\|_2^2=\E_\e\left\|\sum_{A\subseteq
\{1,\ldots,m\}}g_AW_A(\e) S_g(L_g(x_A))\right\|_2^2=\sum_{A\subseteq
\{1,\ldots,m\}} g_A^2\left\|S_g(L_g(x_A))\right\|_2^2,
\end{eqnarray}
and
\begin{eqnarray}\label{eq:orthogonal'}
\E_\e
\left\|S_g'\left(L_g\left(\p_g(\e)\right)\right)\right\|_2^2=\E_\e\left\|\sum_{A\subseteq
\{1,\ldots,m\}}g_AW_A(\e) S_g(L_g(y_A))\right\|_2^2=\sum_{A\subseteq
\{1,\ldots,m\}} g_A^2\left\|S_g(L_g(y_A))\right\|_2^2.
\end{eqnarray}
A combination of the bounds in~\eqref{eq:bilip}
and~\eqref{eq:euclidean} shows that for all $A\subseteq
\{1,\ldots,m\}$ we have $\left\|S_g(L_g(x_A))\right\|_2\le
2D\Delta\left(4K\log(n+1)\right)\|x_A\|$ and for all $\e\in
\{-1,1\}^m$ we have
$\left\|S_g\left(L_g\left(\f_g(\e)\right)\right)\right\|_2\ge
\left\|\f_g(\e)\right\|$. Thus~\eqref{eq:orthogonal} implies that
\begin{eqnarray}\label{eq:before gaussian}
\E_\e \left\|\f_g(\e)\right\|^2\le
4D^2\Delta\left(4K\log(n+1)\right)^2\sum_{A\subseteq
\{1,\ldots,m\}}g_A^2\|x_A\|^2.
\end{eqnarray}
Arguing similarly, while using~\eqref{eq:bilip'},
\eqref{eq:euclidean'} and~\eqref{eq:orthogonal'}, we see that
\begin{eqnarray}\label{eq:before gaussian'}
\E_\e \left\|\p_g(\e)\right\|^2\ge
\frac{1}{4D^2\Delta\left(4K\log(n+1)\right)^2}\sum_{A\subseteq
\{1,\ldots,m\}}g_A^2\|y_A\|^2.
\end{eqnarray}
Taking expectation with respect to the Gaussians
$\{g_A\}_{A\subseteq \{1,\ldots,m\}}$ in~\eqref{eq:before gaussian}
we see that
\begin{eqnarray}\label{eq:before maximizer}
4D^2\Delta\left(4K\log(n+1)\right)^2\sum_{A\subseteq
\{1,\ldots,m\}}\|x_A\|^2\ge
\E_g\E_\e\left\|\sum_{A\subseteq\{1,\ldots,m\}}
g_AW_A(\e)x_A\right\|^2=\E_g\left\|\sum_{A\subseteq\{1,\ldots,m\}}
g_Ax_A\right\|^2,
\end{eqnarray}
where we used the fact that for each fixed $\e\in \{-1,1\}^m$ the
random variables $\{W_A(\e)g_A\}_{A\subseteq \{1,\ldots,m\}}$ have
the same joint distribution as the random variables
$\{g_A\}_{A\subseteq \{1,\ldots,m\}}$. Similarly, taking expectation
in~\eqref{eq:before gaussian'} yields
\begin{eqnarray}\label{eq:before maximizer'}
\sum_{A\subseteq \{1,\ldots,m\}}\|y_A\|^2\le
4D^2\Delta\left(4K\log(n+1)\right)^2\E_g\left\|\sum_{A\subseteq\{1,\ldots,m\}}
g_Ay_A\right\|^2.
\end{eqnarray}
Combining~\eqref{eq:before maximizer} with~\eqref{eq:relabeledx}
and~\eqref{eq:before maximizer'} with~\eqref{eq:relabeledy} we get
the bounds:
$$
T_2(F),C_2(F)\le 2D\Delta\left(4K\log(n+1)\right).
$$
In combination with Kwapien's theorem~\eqref{eq:kwapien} we deduce
that
$$
c_2(F)\le T_2(F)C_2(F)\le 4D^2\Delta\left(4K\log(n+1)\right)^2.
$$
Since $F$ was an arbitrary subspace of $X$ of dimension at most $n$,
it follows that
\begin{eqnarray}\label{eq:recursion}
\Delta(n)\le 4D^2\Delta\left(4K\log(n+1)\right)^2.
\end{eqnarray}
Iterating~\eqref{eq:recursion} $\log^*(n)$ times implies that
$$
\Delta(n)\le 2^{2^{c(K,D)\log^*(n)}},
$$
as required.\qed

\section{Proof of Theorem~\ref{thm:lower}}\label{sec:lower}

%%%%BIGCHANGE
We shall now explain why the $2$-convexification of the Tsirelson
space $T^{(2)}$, as defined in the introduction, satisfies the
%%%%CHANGE
J-L lemma. First we give a definition. Given an increasing sequence
$h(n)$ with $0\le h(n)\le n$, say that a Banach space is
$h$-Hilbertian provided that for every finite dimensional subspace
$E$ of $X$ there are subspaces $F$ and $G$ of $E$ such that
$E=F\oplus G$, $\dim(F)=O(h(\dim E))$ and
$c_2(G)=O(1)$\footnote{Here the direct sum notation means as usual
that for every $y\oplus z\in F\oplus G$ we have $\|y\oplus
z\|_{X}=\Theta\left(\|y\|_{X}+\|z\|_{X}\right)$, where the implied
constants are independent of $E$}.
%This will follow from the following
%fact: If $E\subseteq T^{(2)}$ is an $n$-dimensional subspace then
%there are two subspaces $F,G\subseteq T^{(2)}$ such that $E=F\oplus
%G$, $\dim(F)=O(\log n)$ and $c_2(G)=O(1)$\footnote{Here the direct
%sum notation means as usual that for every $y\oplus z\in F\oplus G$
%we have $\|y\oplus
%z\|_{T^{(2)}}=\Theta\left(\|y\|_{T^{(2)}}+\|z\|_{T^{(2)}}\right)$.}.
%Indeed, assuming the validity of this fact
If the Banach space $X$ is $\log$-Hilbertian, then $X$ satisfies the J-L lemma.  Indeed,
 take $x_1,\ldots,x_n\in
X$ and let $E$ be their span. Write $E=F\oplus G$ as above and
decompose each of the $x_i$ accordingly, i.e. $x_i=y_i\oplus z_i$
where $y_i\in F$ and $z_i\in G$. Since $c_2(G)=O(1)$, by the
J-L lemma we can find a linear operator $L:G\to
G'$, where $G'\subseteq G$ is a subspace of dimension $O(\log n)$,
such that $\|z_i-z_j\|=\Theta(\|L(z_i)-L(z_j)\|)$ for all $i,j\in
\{1,\ldots,n\}$. The linear operator $L':E\to F\oplus G'$ given by
$L'(y\oplus z)=y\oplus L(z)$ has rank $O(\log n)$ and satisfies
$\|x_i-x_j\|=\Theta(\|L'(x_i)-L'(x_j)\|)$, as required.

We now explain  why $T^{(2)}$ satisfies the J-L lemma.
%%%%ENDBIGCHANGE
In~\cite{J76} Johnson defined the following modification of the
Tsirelson space. As in the case of the Tsirelson space, the
construction consists of an inductive definition of a sequence of
norms on $c_{00}$. Once again we set $|||x|||_0=\|x\|_{c_0}$ and
\begin{eqnarray}\label{eq:modifiedtsirl recurse}
|||x|||_{m+1}=\max\left\{
|||x|||_{m},\frac12\sup\left\{\sum_{j=1}^{(n+1)^n}|||P_{A_j}(x)|||_m:\
n\in \N,\ A_1,\ldots,A_{(n+1)^n}\subseteq [n,\infty)\
\mathrm{finite\ \& \ disjoint}\right\}\right\}
\end{eqnarray}
We then define $\|x\|_{\mathscr{T}}\coloneqq \lim_{m\to \infty}
|||x|||_m$, and the modified space $\mathscr{T}^{(2)}$ as the
completion of $c_{00}$ under the norm:
\begin{eqnarray}\label{eq:modifiedT2}
\left\|\sum_{j=1}^\infty x_j
e_j\right\|_{\mathscr{T}^{(2)}}\coloneqq \left\|\sum_{j=1}^\infty
|x_j|^2 e_j\right\|_{\mathscr{T}}^{1/2}.
\end{eqnarray}

%%%%CHANGE
In~\cite{J80} Johnson proved that a certain subspace $Y$
of $\mathscr{T}^{(2)}$ (spanned by a subsequence of the unit vector
basis)
%every $n$-dimensional subspace $E\subseteq Y$ can be
%decomposed as $E=F\oplus G$ for some
%subspaces $F,G\subseteq Y$
%satisfying $\dim(F)=O(\log n)$ and $c_2(G)=O(1)$.
is $\log$-Hilbertian.
In~\cite{CJT84}
Casazza, Johnson and Tzafriri showed that it is not necessary to
pass to the subspace $Y$, and in fact $\mathscr{T}^{(2)}$ itself has
the desired decomposition property. Finally, a deep result of
Casazza and Odell~\cite{CO83} shows that $T^{(2)}$ is just
$\mathscr{T}^{(2)}$ with an equivalent norm. This concludes the
proof of the fact that $T^{(2)}$ satisfies the J-L
lemma.
%%%%ENDCHANGE

It remains to establish the lower bound~\eqref{eq:bellenot}. Note
that the fact that $c_2\left(T^{(2)}\right)=\infty$ already follows
from the original paper of Figiel and Johnson~\cite{FJ74}---our goal
here is to give a quantitative estimate. This will be a simple
consequence of a paper of Bellenot~\cite{Bell84}. Define inductively
a sequence of functions $\{g_k:\N\to \N\}_{k=0}^\infty$ as follows:
$g_0(n)=n+1$ and $g_{i+1}(n)=g_i^{(n)}(n)$, where $g_i^{(n)}$
denotes the $n$-fold iterate of $g_i$, i.e.
$g_{i}^{(n+1)}(j)=g_i(g_i^{(n)}(j))$. The inverse Ackermann function
is the inverse of the function $n\mapsto g_n(n)$, i.e. its value on
$n\in \N$ is the unique integer $k$ such that $g_k(k)<n\le
g_{k+1}(k+1)$. Note that in the literature there are several
variants of the inverse Ackermann function, but it is possible to
show that they are all the same up to bounded additive
%%%%CHANGE
terms---see,
for example,~\cite[Appendix B]{AKNSS08} for a discussion of such
issues. In particular, we define $\alpha(n)$ to be the inverse
of the function $h(n)=g_n(2)$, but its asymptotic behavior is the
same as  the inverse Ackermann function (since $g_n(2)>n$, and therefore
$g_n(n)<g_n(g_n(2))=g_{n+1}(2)$). Now,
%%%%CHANGE
by ~\cite[Proposition 5]{Bell84} for every $k\ge 1$ there exist scalars
$\{x_j\}_{j=1}^{g_{k}(2)}\subseteq \R$ which are not all equal to
$0$ such that
\begin{eqnarray}\label{eq:bad x}
\left\|\sum_{j=1}^{g_{k}(2)}x_j e_j\right\|_T\le
\frac{k}{2^{k}}\sum_{j=3}^{g_{k}(2)}|x_j|.
\end{eqnarray}
Hence, by the definition~\eqref{eq:T2} we have for all
$\e=(\e_1,\ldots,\e_{g_{k}(2)})\in \{-1,1\}^{g_{k}(2)}$,
\begin{eqnarray}\label{eq:bad2}
\left\|\sum_{j=1}^{g_{k}(2)}\e_jx_j e_j\right\|_{T^{(2)}}^2\le
\frac{k^2}{2^{2k}}\sum_{j=1}^{g_{k}(2)}x_j^2.
\end{eqnarray}
Let $F\subseteq T^{(2)}$ denote the span of
$\{e_1,\ldots,e_{g_{k}(2)}\}$. Averaging~\eqref{eq:bad2} over $\e$
and using the definition of the cotype $2$ constant of $F$, we see
that $ C_2(F)\ge 2^k/k$, and therefore the Euclidean distortion of
$F$ is at least $2^k/k$. Since the dimension of $F$ is $g_2(k)$,
this concludes the proof of~\eqref{eq:bellenot}, and hence also the
proof of Theorem~\ref{thm:lower}.\qed

\section{Remarks and open problems}\label{sec:open}

We end this note with some concluding comments and questions that
arise naturally from our work.

\begin{enumerate}
%%%%CHANGE
\item The space
$T^{(2)}$ was the first example of what Pisier~\cite[Chapter 12]{P89}
calls weak Hilbert spaces.  One of the many equivalents for a Banach space $X$
to be a weak Hilbert is that every finite dimensional subspace $E$ of $X$
can be written as $E=F\oplus G$ with $\dim G \ge \delta \dim E$ for some
universal constant $\delta > 0$ and $c_2(G)=O(1)$.  It is not known whether
every weak Hilbert space is $\log $-Hilbertian or even $h$-Hilbertian for some $h(n)=o(n)$.
However, Nielsen and Tomczak-Jaegermann~\cite{NT-J92},
using the same kind of reasoning that works for $T^{(2)}$ (see \cite{CJT84}),
proved that a weak Hilbert space with an unconditional basis is even $2^{O(\alpha(\cdot))}$-Hilbertian.

\item
%{eq:def Delta}
A Banach space $X$ is called asymptotically Hilbertian provied that for each $n$ there is a
finite codimensional subspace $Y$ of $X$ so that $\Delta_Y(n) =O(1)$ ($\Delta_Y(n)$ is defined in
(\ref {eq:def Delta})). Every weak Hilbert space is asymptotically Hilbertian \cite[Chapter 14]{P89}.
The results in \cite{J80} and the argument at the beginning of section \ref{sec:lower} show that every asymptotically Hilbertian space has a subspace which satisfies the J-L lemma.

\item Does there exist a function $f(n )\uparrow \infty$ so that if $X$ is a Banach space for which
%%%%CHANGE
$\Delta(n)=O(f(n))$, where $\Delta(n)$ is as in~\eqref{eq:def
Delta}, i.e. $c_2(E)=O(f(\dim E))$ for all finite dimensional
subspaces $E$ of $X$, then $X$ satisfies the J-L lemma? An
affirmative answer would show that there are natural Banach spaces
other than Hilbert spaces, even some Orlicz sequence spaces, which
satisfy the J-L lemma.

\item A question which obviously arises from our results is to
determine the true rate of ``closeness"
%%%%CHANGE{eq:log*}
(in the sense of (\ref{eq:log*}))
between spaces satisfying
the
J-L lemma and Hilbert space. Which of the
bounds $\Delta(n) = 2^{2^{O(\log^*(n))}}$ and $\Delta(n) =2^{\Omega(\alpha(n))}$ is closer
to the truth?

\item Our argument also  works when the dimension is only assumed to be reduced to a power of
$\log n$, and we get nontrivial bounds even when this dimension is,
say, $2^{(\log n)^{\beta}}$ for some $\beta<1$.
%%%%CHANGE SHOULD MORE BE SAID?
However, except for spaces that
are of type 2 or of cotype 2,
our proof does
not yield any meaningful result when the dimension is lowered to
$n^\gamma$ for some $\gamma\in (0,1)$. The problem is that in the
recursive inequality~\eqref{eq:recursion} the term
$\Delta\left(4K\log(n+1)\right)$ is squared. This happens since in
Kwapien's theorem~\eqref{eq:kwapien} the Euclidean distortion is
bounded by the product of the type $2$ and cotype $2$ constants
rather than by their maximum. While it is tempting to believe that
the true bound in Kwapien's theorem should be
$c_2(X)=O\left(\max\{T_2(X),C_2(X)\}\right)$, it was shown by
%%%%CHANGE
Tomczak-Jaegermann
~\cite[Proposition 2]{T-J83}  that up
to universal constants Kwapien's bound $c_2(X)\le T_2(X)C_2(X)$
cannot be improved.

\item In~\cite{Pisier75} Pisier proved that if a Banach space $X$
satisfies
%%%%CHANGE
 $\Delta(n)=o(\log n)$, then $X$ is superreflexive;  i.e.,  $X$ admits an
equivalent norm which is uniformly convex. Hence any space
satisfying the assumptions of Theorem~\ref{thm:upper} is
superreflexive.

\item It is of interest to study dimension reduction into arbitrary
low dimensional normed spaces, since this can serve just as well for
the purpose of data compression (see~\cite{Mat90}). Assume that $X$
is a Banach space such that for every $n$ and $x_1,\ldots,x_n\in X$
there exists a $d(n)$-dimensional Banach space $Y$ and a linear
mapping $L:X\to Y$ such that $\|x_i-x_j\|\le\|L(x_i)-L(x_j)\|\le
D\|x_i-x_j\|$ for all $i,j\in \{1,\ldots,n\}$. Since by John's
theorem~\cite{John48} we have $c_2(Y)\le \sqrt{d(n)}$ we can argue
similarly to the proof in Section~\ref{sec:upper} (in this simpler
case the proof is close to the argument in~\cite{LMN05}), while
using the result of Tomczak-Jaegermann~\eqref{eq:T-J}, to deduce
that $T_2(F),C_2(F)\le 2\pi\sqrt{d(n)}$. By Kwapien's theorem we
deduce that $c_2(F)\le 4\pi^2 d(n)$. If $d(n)\le n^{\gamma}$ for
some $\gamma\in (0,1)$ which is independent of $n$ and $F$, the fact
that $T_2(F),C_2(F)\le 2\pi n^{\gamma/2}$ for every $n$-dimensional
subspace $F\subseteq X$ implies (see~\cite{T-J89}) that $X$ has type
$\frac{2}{1+\gamma}-\e$ and cotype $\frac{2}{1-\gamma}+\e$ for every
$\e>0$. In particular, if $d(n)=n^{o(1)}$ then $X$ has type $2-\e$
and cotype $2+\e$ for every $\e>0$.

\item We do not know of any non-trivial linear dimension reduction
result in $L_p$ for $p\in [1,\infty)\setminus\{2\}$. For example, is
it possible to embed with $O(1)$ distortion via a linear mapping any
$n$-point subset of $L_1$ into a subspace of $L_1$ of dimension,
say, $n/4$, or even into $\ell_1^{n/4}$? Remarkably even such modest
goals seem to be beyond the reach of current techniques. Clearly
$n$-point subsets of $L_1$ are in their $n$-dimensional span, but we
do not know if they embed with constant distortion into $\ell_1^d$
when $d=O(n)$. Schechtman proved in~\cite{Schecht87} that we can
take $d=O(n\log n)$. We refer to~\cite{Schecht87,BLM89,Tal90,JS01}
for more information on the harder problem of embedding
$n$-dimensional subspaces of $L_1$ into low dimensional $\ell_1^d$.
We also refer to these references for similar results in $L_p$
spaces.

\end{enumerate}

%\cite{NT-J92,P89}

%\bibliographystyle{abbrv}
%%%%CHANGE
%\bibliography{JL-L3}

\begin{thebibliography}{10}

\bibitem{AC06}
N.~Ailon and B.~Chazelle.
\newblock Approximate nearest neighbors and the fast {J}ohnson-{L}indenstrauss
  transform.
\newblock In {\em STOC 2006: ACM Symposium on Theory of Computing}, pages
  557--563. 2006.

\bibitem{AKNSS08}
N.~Alon, H.~Kaplan, G.~Nivasch, M.~Sharir, and S.~Smorodinsky.
\newblock Weak {$\epsilon$}-nets and interval chains.
\newblock In {\em Proceedings of the Nineteenth Annual {ACM}-{SIAM} Symposium
  on Discrete Algorithms, {SODA} 2008, San Francisco, California, {USA},
  January 20-22, 2008}, pages 1194--1203. SIAM, 2008.

\bibitem{AR92}
J.~Arias-de Reyna and L.~Rodr{\'{\i}}guez-Piazza.
\newblock Finite metric spaces needing high dimension for {L}ipschitz
  embeddings in {B}anach spaces.
\newblock {\em Israel J. Math.}, 79(1):103--111, 1992.

\bibitem{Ball97}
K.~Ball.
\newblock An elementary introduction to modern convex geometry.
\newblock In {\em Flavors of geometry}, volume~31 of {\em Math. Sci. Res. Inst.
  Publ.}, pages 1--58. Cambridge Univ. Press, Cambridge, 1997.

\bibitem{Bell84}
S.~F. Bellenot.
\newblock The {B}anach space {$T$} and the fast growing hierarchy from logic.
\newblock {\em Israel J. Math.}, 47(4):305--313, 1984.

\bibitem{BL00}
Y.~Benyamini and J.~Lindenstrauss.
\newblock {\em Geometric nonlinear functional analysis. {V}ol. 1}, volume~48 of
  {\em American Mathematical Society Colloquium Publications}.
\newblock American Mathematical Society, Providence, RI, 2000.

\bibitem{Bou85}
J.~Bourgain.
\newblock  On Lipschitz embedding of finite metric spaces in Hilbert space.
\newblock {\em Israel J. Math.},  52(1-2):46--52, 1985.

\bibitem{BLM89}
J.~Bourgain, J.~Lindenstrauss, and V.~Milman.
\newblock Approximation of zonoids by zonotopes.
\newblock {\em Acta Math.}, 162(1-2):73--141, 1989.

\bibitem{BC05}
B.~Brinkman and M.~Charikar.
\newblock On the impossibility of dimension reduction in {$l\sb 1$}.
\newblock {\em J. ACM}, 52(5):766--788 (electronic), 2005.

\bibitem{CJT84}
P.~G. Casazza, W.~B. Johnson, and L.~Tzafriri.
\newblock On {T}sirelson's space.
\newblock {\em Israel J. Math.}, 47(2-3):81--98, 1984.

\bibitem{CO83}
P.~G. Casazza and E.~Odell.
\newblock Tsirelson's space and minimal subspaces.
\newblock In {\em Texas functional analysis seminar 1982--1983 (Austin, Tex.)},
  Longhorn Notes, pages 61--72. Univ. Texas Press, Austin, TX, 1983.

\bibitem{CS89}
P.~G. Casazza and T.~J. Shura.
\newblock {\em Tsirel\cprime son's space}, volume 1363 of {\em Lecture Notes in
  Mathematics}.
\newblock Springer-Verlag, Berlin, 1989.
\newblock With an appendix by J. Baker, O. Slotterbeck and R. Aron.

\bibitem{CS02}
M.~Charikar and A.~Sahai.
\newblock Dimension reduction in the $\ell_1$ norm.
\newblock In {\em 43rd Annual IEEE Conference on Foundations of Computer
  Science}, pages 251--260. IEEE Computer Society, 2002.

\bibitem{Ts74}
B.~S. Cirel{\cprime}son.
\newblock It is impossible to imbed {$1\sb{p}$} of {$c\sb{0}$} into an
  arbitrary {B}anach space.
\newblock {\em Funkcional. Anal. i Prilo\v zen.}, 8(2):57--60, 1974.

\bibitem{DGK63}
L.~Danzer, B.~Gr{\"u}nbaum, and V.~Klee.
\newblock Helly's theorem and its relatives.
\newblock In {\em Proc. Sympos. Pure Math., Vol. VII}, pages 101--180. Amer.
  Math. Soc., Providence, R.I., 1963.

\bibitem{FJ74}
T.~Figiel and W.~B. Johnson.
\newblock A uniformly convex {B}anach space which contains no {$l\sb{p}$}.
\newblock {\em Compositio Math.}, 29:179--190, 1974.

\bibitem{FLM77}
T.~Figiel, J.~Lindenstrauss, and V.~D. Milman.
\newblock The dimension of almost spherical sections of convex bodies.
\newblock {\em Acta Math.}, 139(1-2):53--94, 1977.

\bibitem{Indyk01}
P.~Indyk.
\newblock Algorithmic applications of low-distortion geometric embeddings.
\newblock In {\em 42nd IEEE Symposium on Foundations of Computer Science (Las
  Vegas, NV, 2001)}, pages 10--33. IEEE Computer Soc., Los Alamitos, CA, 2001.

\bibitem{Indyk06}
P.~Indyk.
\newblock Stable distributions, pseudorandom generators, embeddings, and data
  stream computation.
\newblock {\em J. ACM}, 53(3):307--323 (electronic), 2006.

\bibitem{IM99}
P.~Indyk and R.~Motwani.
\newblock Approximate nearest neighbors: towards removing the curse of
  dimensionality.
\newblock In {\em STOC '98 (Dallas, TX)}, pages 604--613. ACM, New York, 1999.

\bibitem{John48}
F.~John.
\newblock Extremum problems with inequalities as subsidiary conditions.
\newblock In {\em Studies and Essays Presented to R. Courant on his 60th
  Birthday, January 8, 1948}, pages 187--204. Interscience Publishers, Inc.,
  New York, N. Y., 1948.

\bibitem{J76}
W.~B. Johnson.
\newblock A reflexive {B}anach space which is not sufficiently {E}uclidean.
\newblock {\em Studia Math.}, 55(2):201--205, 1976.

\bibitem{J80}
W.~B. Johnson.
\newblock Banach spaces all of whose subspaces have the approximation property.
\newblock In {\em Special topics of applied mathematics (Proc. Sem., Ges. Math.
  Datenverarb., Bonn, 1979)}, pages 15--26. North-Holland, Amsterdam, 1980.

\bibitem{JL84}
W.~B. Johnson and J.~Lindenstrauss.
\newblock Extensions of {L}ipschitz mappings into a {H}ilbert space.
\newblock In {\em Conference in modern analysis and probability (New Haven,
  Conn., 1982)}, volume~26 of {\em Contemp. Math.}, pages 189--206. Amer. Math.
  Soc., Providence, RI, 1984.

  \bibitem{JLS87}
W.~B. Johnson, J.~Lindenstrauss, and G.~Schechtman.
\newblock On Lipschitz embedding of finite metric spaces In low-dimensional normed spaces.
 \newblock In {\em Geometrical aspects of functional analysis (1985/86)},
 pages 177--184.
  \newblock Lecture Notes in Math., 1267, Springer, Berlin, 1987.

\bibitem{JS01}
W.~B. Johnson and G.~Schechtman.
\newblock Finite dimensional subspaces of {$L\sb p$}.
\newblock In {\em Handbook of the geometry of Banach spaces, Vol. I}, pages
  837--870. North-Holland, Amsterdam, 2001.

\bibitem{Kleinberg97}
J.~M. Kleinberg.
\newblock Two algorithms for nearest-neighbor search in high dimensions.
\newblock In {\em STOC '97 (El Paso, TX)}, pages 599--608 (electronic). ACM,
  New York, 1999.

\bibitem{KOR00}
E.~Kushilevitz, R.~Ostrovsky, and Y.~Rabani.
\newblock Efficient search for approximate nearest neighbor in high dimensional
  spaces.
\newblock {\em SIAM J. Comput.}, 30(2):457--474 (electronic), 2000.

\bibitem{Kw72}
S.~Kwapie{\'n}.
\newblock Isomorphic characterizations of inner product spaces by orthogonal
  series with vector valued coefficients.
\newblock {\em Studia Math.}, 44:583--595, 1972.
\newblock Collection of articles honoring the completion by Antoni Zygmund of
  50 years of scientific activity, VI.

\bibitem{LMN05}
J.~R. Lee, M.~Mendel, and A.~Naor.
\newblock Metric structures in {$L\sb 1$}: dimension, snowflakes, and average
  distortion.
\newblock {\em European J. Combin.}, 26(8):1180--1190, 2005.

\bibitem{LN04}
J.~R. Lee and A.~Naor.
\newblock Embedding the diamond graph in {$L\sb p$} and dimension reduction in
  {$L\sb 1$}.
\newblock {\em Geom. Funct. Anal.}, 14(4):745--747, 2004.

\bibitem{Mat90}
J.~Matou{\v{s}}ek.
\newblock Bi-{L}ipschitz embeddings into low-dimensional {E}uclidean spaces.
\newblock {\em Comment. Math. Univ. Carolin.}, 31(3):589--600, 1990.

\bibitem{Mat96}
J.~Matou{\v{s}}ek.
\newblock On the distortion required for embedding finite metric spaces into
  normed spaces.
\newblock {\em Israel J. Math.}, 93:333--344, 1996.

\bibitem{NT-J92}
N.~J. Nielsen and N.~Tomczak-Jaegermann.
\newblock Banach lattices with property ({H}) and weak {H}ilbert spaces.
\newblock {\em Illinois J. Math.}, 36(3):345--371, 1992.

\bibitem{Pisier75}
G.~Pisier.
\newblock Martingales with values in uniformly convex spaces.
\newblock {\em Israel J. Math.}, 20(3-4):326--350, 1975.

\bibitem{Pisier86}
G.~Pisier.
\newblock {\em Factorization of linear operators and geometry of {B}anach
  spaces}, volume~60 of {\em CBMS Regional Conference Series in Mathematics}.
\newblock Published for the Conference Board of the Mathematical Sciences,
  Washington, DC, 1986.

\bibitem{P89}
G.~Pisier.
\newblock {\em The volume of convex bodies and {B}anach space geometry},
  volume~94 of {\em Cambridge Tracts in Mathematics}.
\newblock Cambridge University Press, Cambridge, 1989.

\bibitem{Schecht87}
G.~Schechtman.
\newblock More on embedding subspaces of {$L\sb p$} in {$l\sp n\sb r$}.
\newblock {\em Compositio Math.}, 61(2):159--169, 1987.

\bibitem{Tal90}
M.~Talagrand.
\newblock Embedding subspaces of {$L\sb 1$} into {$l\sp N\sb 1$}.
\newblock {\em Proc. Amer. Math. Soc.}, 108(2):363--369, 1990.

\bibitem{T-J79}
N.~Tomczak-Jaegermann.
\newblock Computing {$2$}-summing norm with few vectors.
\newblock {\em Ark. Mat.}, 17(2):273--277, 1979.

\bibitem{T-J83}
N.~Tomczak-Jaegermann.
\newblock The {B}anach-{M}azur distance between symmetric spaces.
\newblock {\em Israel J. Math.}, 46(1-2):40--66, 1983.

\bibitem{T-J89}
N.~Tomczak-Jaegermann.
\newblock {\em Banach-{M}azur distances and finite-dimensional operator
  ideals}, volume~38 of {\em Pitman Monographs and Surveys in Pure and Applied
  Mathematics}.
\newblock Longman Scientific \& Technical, Harlow, 1989.

\bibitem{Vem04}
S.~S. Vempala.
\newblock {\em The random projection method}.
\newblock DIMACS Series in Discrete Mathematics and Theoretical Computer
  Science, 65. American Mathematical Society, Providence, RI, 2004.
\newblock With a foreword by Christos H.\ Papadimitriou.

\end{thebibliography}

\def\cprime{$'$} \def\cprime{$'$} \def\cprime{$'$}

\end{document}